\documentclass[11pt]{amsart}

\usepackage{amssymb, latexsym}
\theoremstyle{plain}
\newtheorem*{theorem}{Theorem}

\theoremstyle{remark}

\newtheorem{remark}{Remark}

\numberwithin{equation}{section}

\begin{document}

\title {Uniqueness/nonuniqueness for nonnegative solutions
of the Cauchy problem for $u_t=\Delta u-u^p$ in a  punctured  space}

\author{Ross G. Pinsky}
\address{Department of Mathematics\\
Technion---Israel Institute of Technology\\
Haifa, 32000\\
Israel}
\email{pinsky@math.technion.ac.il}

\subjclass[2000]{Primary; 60J80, 60J60}
\date{}
\begin{abstract}

Consider classical solutions to the following Cauchy problem  in a punctured space:
\begin{equation}\label{abstract}
\begin{aligned}
&u_t=\Delta u
-u^p\ \ \text{in}\  (R^n-\{0\})\times(0,\infty);\\
& u(x,0)=g(x)\ge0 \ \ \text{in}\ R^n-\{0\};\\
&u\ge0 \ \ \text{in} \ (R^n-\{0\})\times[0,\infty).\\
\end{aligned}
\end{equation}
We prove that if $p\ge\frac n{n-2}$, then
the solution to \eqref{abstract} is unique
for each $g$. On the other hand, if $p<\frac n{n-2}$, then
uniqueness does not hold when $g=0$; that is, there exists
a nontrivial solution with vanishing initial data.

\end{abstract}
\maketitle
\section{Introduction and Statement of Results}\label{S:intro}

The study of uniqueness
 in the class of all classical solutions for the Cauchy problem
\begin{equation}\label{gen}
\begin{aligned}
&u_t=Lu+f(x,u)\ \ \text{in}\ R^n\times(0,\infty);\\
&u(x,0)=g(x)\ge0\  \ \text{in}\ R^n;\\
&u\ge0\ \ \text{in} \ R^n\times[0,\infty),
\end{aligned}
\end{equation}
where $L$ is a second order elliptic operator,
goes back to Brezis \cite{B84}, where uniqueness was proved
in the case that $L=\Delta$ and $f(x,u)=-u^p$ with $p>1$.
In recent papers \cite{EP03, P04}, the dichotomy
between uniqueness/nonuniqueness was investigated
for general second order elliptic operators
$L$ and quite general nonlinearities $f$, which approach $-\infty$ at a superlinear
rate as $u\to\infty$.
We emphasize that  uniqueness here  is with regard to the   class
of \it all\rm\ classical solutions
 to \eqref{gen} (with no growth restrictions).
For example, let  $L=\Delta$ and
$f(x,u)=-\gamma(x)u^p$.
In \cite{EP03} it was proved that
uniqueness holds for \eqref{gen} if
 $\gamma(x)\ge c_1\exp(-c_2|x|^2)$,
for some $c_1,c_2>0$, while uniqueness does not hold
in \eqref{gen} with initial data $g=0$ if
 $\gamma(x)\le c\exp(-|x|^{2+\epsilon})$,
 for some $c>0$. On the other hand, if one  looks only at \it mild\rm\
 solutions, then it is well known that uniqueness holds above for all bounded
 $\gamma$ \cite{P83}.

In this paper, we study the question of uniqueness
for the same semilinear equation $u_t=\Delta u-u^p$ studied by Brezis, but
replace the space $R^n$  by
the \it punctured space\rm\
$R^n-\{0\}$, $n\ge2$, thus allowing for unboundedness of solutions
in a neighborhood of 0  at all times $t\ge0$:
\begin{equation}\label{semilinear}
\begin{aligned}
&u_t=\Delta u-u^p\  \ \text{in}\ (R^n-\{0\})\times(0,\infty);\\
& u(x,0)=g(x)\ge0 \  \ \text{in}\ R^n-\{0\};\\
&u\ge0 \ \ \text{in}\ (R^n-\{0\})\times[0,\infty).\\
\end{aligned}
\end{equation}
We assume that $g\in C(R^n-\{0\})$.

We prove the following theorem.
\begin{theorem}
\
\begin{enumerate}
\item Let $p<\frac n{n-2}$. Then
there exists a nontrivial solution to \eqref{semilinear} with
initial data $g=0$.
\item Let $p\ge\frac n{n-2}$. Then
for each $g$ there exists a unique solution to \eqref{semilinear}.
\end{enumerate}
\end{theorem}
\begin{remark}
For the case $p\in(1,2]$, this  result has an important
interpretation  with regard to the theory of super-Brownian motion;
see \cite{P04a} for details.

\end{remark}

\begin{remark}
Brezis and Friedman \cite{BF83}  studied the problem
$u_t=\Delta u-u|u|^{p-1}$ in $R^n\times(0,\infty)$, with
the initial condition $u(x)=\delta_0(x)$, the Dirac $\delta$-function
at 0. They showed that a solution exists if and only
if $p<\frac{n+2}n$. More recently, Marcus and Veron \cite{MV99} have shown that
for positive solutions of the above equation, the same
critical exponent appears when one allows for even more
singular initial conditions---namely, not necessarily
locally bounded Borel measures.
In these papers, the solution is required to be classical at
$x=0$, for $t>0$, whereas the present paper deals with
the situation in which $x=0$ is excluded for all times $t\ge0$.
As such, it is `easier'' to obtain nontrivial solutions
in the present case, and this is
manifested through the larger critical exponent, $\frac n{n-2}$
as compared to $\frac{n+2}n$.
\end{remark}

\begin{remark}
Note that for $n=2$, nonuniqueness prevails for the problem
in this paper  for all $p>1$.

\end{remark}

We give a very simple proof of part (1) of the theorem by exploiting
a recent result in \cite{P04}. For the proof of part (2), we construct appropriate
supersolutions.

\section{Proof of Theorem}
We begin by noting that existence follows by standard methods;
see \cite{L96} or  \cite{EP03} (this latter reference treats the
case that the domain is $R^n$, but the same techniques work
in the punctured space). Thus, it remains to consider uniqueness.
\begin{proof}[Proof of part (1).]
Since the problem
is radially symmetric, it suffices to show that uniqueness fails
for the radially symmetric equation
\begin{equation}\label{onedim}
\begin{aligned}
&u_t=u_{rr}+\frac {n-1}ru_r-u^p,\ r\in(0,\infty), \ t>0;\\
&
u(r,0)=0, \ r\in(0,\infty);\\
&u\ge0,\   r\in(0,\infty), \ t\ge0.
\end{aligned}
\end{equation}
By assumption, we have  $p<\frac n{n-2}$,
or equivalently, $n<\frac{2p}{p-1}$.
Thus, the function $W(x)=Cx^{-\frac2{p-1}}$, where
$C^{p-1}=\frac2{p-1}(\frac{2p}{p-1}-n)$,
is a positive, stationary solution
of the parabolic equation
$u_t=u_{rr}+\frac {n-1}ru_r-u^p$ in $(0,\infty)$.
By \cite[Theorem 2-ii]{P04},
the fact that there exists a nontrivial positive, stationary solution
guarantees that uniqueness does not hold for the corresponding
parabolic equation with initial data 0; that is, uniqueness does not
hold for \eqref{onedim}.
Actually, the result in
\cite{P04} is for equations with domain $R^n$, $n\ge1$, whereas
the domain here is $(0,\infty)$. One can check that the proof also
holds in a half space, but more simply, one can make the change of
variables $z=\frac1x-x$, which converts the problem to all of $R$.

\noindent \it Proof of part (2).\rm\
We will prove uniqueness for \eqref{semilinear} in the case
of vanishing initial data. This is enough because
by \cite[Proposition 3]{EP03}, uniqueness for arbitrary $g\ge0$
follows from uniqueness for the case $g=0$.
We write the condition $p\ge\frac n{n-2}$ in the form
$n\ge\frac{2p}{p-1}$.
For technical reasons, it will be necessary to treat the cases
$n>\frac{2p}{p-1}$ and $n=\frac{2p}{p-1}$ separately.

We first consider the case $n>\frac{2p}{p-1}$.
For $\epsilon$ and $R$ satisfying
$0<\epsilon<1$ and $R>1$, and some $l\in(0,1]$, define
\begin{equation}\label{phi}
\phi_{R,\epsilon}(x)=((|x|-\epsilon)(R-|x|))^{-\frac2{p-1}}(1+|x|)^\frac2{p-1}
(1+\frac{\epsilon^l}{ |x|^l} R^\frac2{p-1}).
\end{equation}
Also, for $R$ and $\epsilon$ as above, and some $\gamma>0$, define
\begin{equation}\label{psi}
\psi_{R,\epsilon}(x,t)=\phi_{R,\epsilon}(x)\exp(\gamma(t+1)).
\end{equation}
Note that $\psi_{R,\epsilon}(x,0)>0$, for
$|x|\in(\epsilon, R)$, and  $\psi_{R,\epsilon}(x,t)
=\infty$, for $|x|=\epsilon$ or $|x|=R$.
We will show that for
all sufficiently large $R$ and all sufficiently  small $\epsilon$, and for
$\gamma$ sufficiently large and $l$ sufficiently small, independent of those
 $R$ and $\epsilon$, one has
\begin{equation}\label{supersolution}
\Delta\psi_{R,\epsilon}-\psi^p_{R,\epsilon}
-(\psi_{R,\epsilon})_t\le0,\ \ \text{for} \
\epsilon<|x|<R\ \ \text{and}\ t>0.
\end{equation}
It then follows from the maximum principle for semi-linear equations
that  every  solution $u(x,t)$ to \eqref{semilinear} satisfies
\begin{equation}\label{upsi}
u(x,t)\le\psi_{R,\epsilon}(x,t), \ \ \text{for}\ \epsilon<|x|<R
\ \ \text{and} \ t\in[0,\infty).
\end{equation}
Substituting \eqref{phi} and \eqref{psi} in \eqref{upsi}, letting
$\epsilon\to0$, and then letting $R\to\infty$, one concludes  that
$u(x,t)\equiv0$.
Thus, it remains to show \eqref{supersolution}.

From now on we will use radial coordinates, writing $\phi(r)$
for $\phi(x)$ with $|x|=r$, and similarly for $\psi$.
We have
\begin{equation}\label{first}
\begin{aligned}
&\exp(-\gamma(t+1))(\psi_{R,\epsilon})_r=\\
&-(\frac2{p-1})((r-\epsilon)(R-r))^{-\frac2{p-1}-1}
(R+\epsilon-2r)(1+r)^\frac2{p-1}(1+\frac{\epsilon^l}{r^l}R^\frac2{p-1})\\
&+(\frac2{p-1})((r-\epsilon)(R-r))^{-\frac2{p-1}}(1+r)^{\frac2{p-1}-1}
(1+\frac{\epsilon^l}{r^l}R^\frac2{p-1})\\
&-l((r-\epsilon)(R-r))^{-\frac2{p-1}}(1+r)^{\frac2{p-1}}
\frac{\epsilon^l}{r^{l+1}}R^\frac2{p-1},
\end{aligned}
\end{equation}
and
\begin{equation}\label{second}
\begin{aligned}
&\exp(-\gamma(t+1))\big((r-\epsilon)(R-r)\big)^{-\frac2{p-1}-2}
(\psi_{R,\epsilon})_{rr}=\\
&(\frac2{p-1})(\frac2{p-1}+1)(R+\epsilon-2r)^2(1+r)^\frac2{p-1}
(1+\frac{\epsilon^l}{r^l}R^\frac2{p-1})\\
&+2(\frac2{p-1})(r-\epsilon)
(R-r)(1+r)^{\frac2{p-1}}(1+\frac{\epsilon^l}{r^l}R^\frac2{p-1})\\
&-2(\frac2{p-1})^2(r-\epsilon)(R-r)(R+\epsilon-2r)(1+r)^{\frac2{p-1}-1}
(1+\frac{\epsilon^l}{r^l}R^\frac2{p-1})\\
&+2l(\frac2{p-1})(r-\epsilon)(R-r)(R+\epsilon-2r)(1+r)^{\frac2{p-1}}
\frac{\epsilon^l}{r^{l+1}}R^{\frac2{p-1}}\\
&+(\frac2{p-1})(\frac2{p-1}-1)((r-\epsilon)(R-r))^2
(1+r)^{\frac2{p-1}-2}
(1+\frac{\epsilon^l}{r^l}R^\frac2{p-1})\\
&-2l(\frac2{p-1})((r-\epsilon)(R-r))^2
(1+r)^{\frac2{p-1}-1}
\frac{\epsilon^l} {r^{l+1}}R^\frac2{p-1}\\
&+l(l+1)((r-\epsilon)(R-r))^2(1+r)^\frac2{p-1}
\frac{\epsilon^l}{r^{l+2}}R^\frac2{p-1}.
\end{aligned}
\end{equation}
Using \eqref{phi}, \eqref{psi}, \eqref{first} and the fact that $\frac2{p-1}+2=\frac{2p}{p-1}$,
we have
\begin{equation}\label{other}
\begin{aligned}
&\exp(-\gamma(t+1))\big((r-\epsilon)(R-r)\big)^{-\frac2{p-1}-2}
\Big(\frac {n-1}r(\psi_{R,\epsilon})_r-\psi_{R,\epsilon}^p
-(\psi_{R,\epsilon})_t\Big)=\\
&-(\frac2{p-1})(\frac {n-1}r)(r-\epsilon)(R-r)(R+\epsilon-2r)
(1+r)^{\frac2{p-1}}(1+\frac{\epsilon^l}{r^l}R^\frac2{p-1})\\
&+(\frac2{p-1})(\frac {n-1}r)((r-\epsilon)(R-r))^2
(1+r)^{\frac2{p-1}-1}(1+\frac{\epsilon^l}{r^l}R^\frac2{p-1})\\
&-l(\frac {n-1}r)((r-\epsilon)(R-r))^2
(1+r)^{\frac2{p-1}}\frac{\epsilon^l} {r^{l+1}}R^\frac2{p-1}\\
&-\gamma((r-\epsilon)(R-r))^2(1+r)^{\frac2{p-1}}
(1+\frac{\epsilon^l}{r^l}R^\frac2{p-1})\\
&-(1+r)^{\frac{2p}{p-1}}(1+\frac{\epsilon^l}{r^l}R^\frac2{p-1})^p
\exp((p-1)\gamma(t+1)).
\end{aligned}
\end{equation}
We will   show that for all sufficiently large $R$ and sufficiently small
$\epsilon$, and
for $\gamma$ sufficiently large and $l$ sufficiently small,
independent of those $R$ and
$\epsilon$, the sum of the right hand sides of \eqref{second}
and \eqref{other} is non-positive.
This will  prove \eqref{supersolution}.

We  denote the seven terms on the right hand side of
\eqref{second} by $J_1-J_7$, and the five terms on the right hand
side of \eqref{other} by $I_1-I_5$.
Note that the terms
that are positive are $J_1, J_2, J_4, J_5, J_7$ and $I_2$.
In what follows, $M$ will denote a positive number that can be made
as large as one desires  by choosing $\gamma$ sufficiently large.
Consider first those $r$ satisfying $r\ge cR$, where $c$ is a fixed positive
number. For $r$ in this range, we have
$|I_5|\ge MR^{\frac2{p-1}+2}(1+\epsilon^l R^{\frac2{p-1}-l})$.
It is easy to see that for $M$ sufficiently large,
 $|I_5|$ dominates each of the positive terms,
uniformly over large $R$ and small $\epsilon$,
and thus (since $M$ can be made arbitrarily large) also the sum of
all of the positive terms. Now consider those $r$ for which
$\delta_0\le r\le C$, for some constants $0<\delta_0<C$.
For $r$ in this range and $\epsilon$ sufficiently small, we have
$|I_4|\ge MR^2(1+\epsilon^l R^{\frac2{p-1}})$, and it is easy to see
that for $M$ sufficiently large, $|I_4|$ dominates
each of the positive terms, uniformly over large $R$ and small $\epsilon$,
and thus, also the sum of all of the positive terms.
One can also show that the transition from $r$ of order unity to
$r$ of order $R$ causes no problem. Thus, we conclude that for any fixed
$\delta_0>0$ and $\gamma$ sufficiently large, the sum of the  right hand sides of
\eqref{second} and \eqref{other} is negative for all large $R$
and small $\epsilon$.
Note that all this holds uniformly over $l\in(0,1]$. The parameter $l$
has not been needed yet.

We now turn to the delicate
situation---when $\epsilon\le r\le \delta_0$.
For later use, we remind the reader that $\delta_0$ may be chosen as small
as one likes.
(Note that at $r=\epsilon$, all the terms vanish except
$J_1$ and $I_5$.  Using the fact that
$\frac2{p-1}+2=\frac{2p}{p-1}$,
it is easy to see that for sufficiently large
$\gamma$, $|I_5(\epsilon)|$ dominates $J_1(\epsilon)$, uniformly over all
large $R$ and small $\epsilon$. However, when
$r$ is small, but on an  order larger  than $\epsilon$, the analysis becomes
a lot more involved.)
In the sequel, whenever we say that a condition holds for $\gamma$
or $M$  sufficiently large, or for $l$ sufficiently small,
we mean independent of $R$ and $\epsilon$.

We first take care of the easy terms. Clearly, $J_5\le |I_4|$
if $\gamma$ is sufficiently large. Also $J_7=\frac{l+1}{n-1} |I_3|\le |I_3|$,
if $l$ is chosen sufficiently small.
(This last inequality holds since by assumption, $n>\frac{2p}{p-1}$; thus,
 $n>2$ for all choices of $p$.)

We now show that for $\gamma$ sufficiently
large,  $J_2\le |I_4|+|I_5|$, for $\epsilon\le r\le\delta_0$.
(We are reusing $|I_4|$ here. Later we will reuse $|I_5|$. This is
permissible because $\gamma$
can be chosen as large as we like.)
To show this inequality, it suffices to show that for $M$ sufficiently large,
\begin{equation}\label{J_2}
(r-\epsilon)R
\le M(r-\epsilon)^2R^2
+M(1+\frac{\epsilon^l}{r^l}R^\frac2{p-1})^{p-1},\
\text{for}\ r\in[\epsilon,\delta_0]
\end{equation}
A trivial calculation shows that the left hand side of \eqref{J_2}
is less than the first term on the right hand side if
$r\ge\epsilon+\frac1{RM}$. If $r\in[\epsilon,\epsilon+\frac1{RM}]$,
then the left hand side of \eqref{J_2} is less than
or equal to $\frac1M$ while the second term on the right hand side
is greater than $M$. We conclude that \eqref{J_2}
holds with $M\ge1$.

It remains to consider $J_1$, $J_4$ and $I_2$.
We will show that for $\gamma$ sufficiently large,
\begin{equation}\label{J_1I_1I_5}
J_1+J_4+I_2+I_1+I_5\le0, \ \text{for}\ r\in[\epsilon,\delta_0].
\end{equation}
Since $I_2$ has the factor $(r-\epsilon)^2$,
while $I_1$ has the factor $(r-\epsilon)$,
and since $\frac{R-r}{R+\epsilon-2r}$ can be made
arbitrarily close to 1 by choosing $R$ sufficiently large,
it follows that  for any $\eta>0$,  we have
$I_2\le\eta|I_1|$, for $r\in[\epsilon,\delta_0]$,
if we choose $\delta_0$ sufficiently small
and $R$ sufficiently large.
Note that $J_4\le\frac {2l}{n-1}|I_1|$.
Thus,
\begin{equation}\label{J_1}
J_1+J_4+I_2+I_1\le J_1+(1-\frac{2l}{n-1}-\eta)I_1=J_1+(1-\kappa)I_1,
\end{equation}
where $\kappa=\frac{2l}{n-1}+\eta$. Also note that since we are free
to choose $l$ and $\eta$ as small as we like, the same holds for
$\kappa$.
We have
\begin{equation}\label{J_1I_1}
\begin{aligned}
&J_1+(1-\kappa)I_1=
(1+r)^{\frac2{p-1}}
(R+\epsilon-2r)
(1+\frac{\epsilon^l}{r^l}R^\frac2{p-1})
\times\\
&\Big(\frac{2(p+1)}{(p-1)^2}(R+\epsilon-2r)-(1-\kappa)\frac{2(n-1)}{p-1}
(\frac{r-\epsilon}r)(R-r)\Big).
\end{aligned}
\end{equation}
From the assumption that $n>\frac{2p}{p-1}$, it
follows that for $\kappa$ sufficiently small and $R$ sufficiently large,
\begin{equation}\label{factor}
\Big(\frac{2(p+1)}{(p-1)^2}(R+\epsilon-2r)-(1-\kappa)
\frac{2(n-1)}{p-1}(\frac{r-\epsilon}r)(R-r)\Big)
\le C\frac\epsilon r R, \ r\in[\epsilon,\delta_0],
\end{equation}
for some $C>0$.
From \eqref{J_1}-\eqref{factor}, we obtain
\begin{equation}\label{estimate}
\begin{aligned}&J_1+J_4+I_2+I_1\le
C\frac\epsilon r R(1+r)^{\frac2{p-1}}
(R+\epsilon-2r)(1+\frac{\epsilon^l}{r^l}R^\frac2{p-1}),\\
&\text{for}\ r\in[\epsilon,\delta_0].
\end{aligned}
\end{equation}

In light of \eqref{estimate}, in order to prove \eqref{J_1I_1I_5},
it suffices to show that
\begin{equation}\label{final}
\frac\epsilon r R^2\le M(1+\frac{\epsilon^l}{r^l}R^\frac2{p-1})^{p-1},
\ r\in[\epsilon,\delta_0],
\end{equation}
for sufficiently large $M$.
Choose $l$ sufficiently small so that
 $l(p-1)\le1$.
 Then  the right hand side of \eqref{final}
is greater or equal to $M\frac\epsilon r R^2$.

We now turn to the case $n=\frac{2p}{p-1}$.
For $\epsilon$ and $R$ satisfying
$0<\epsilon<1$ and $R>1$, and some $c\ge2$, define
\begin{equation}\label{phiagain}
\phi_{R,\epsilon}(x)=((|x|-\epsilon)(R-|x|))^{-\frac2{p-1}}(1+|x|)^\frac2{p-1}
\big(1+(\frac{R^2}{\log\frac {c|x|}\epsilon})^\frac1{p-1}\big).
\end{equation}
Note that the only difference between $\phi_{R,\epsilon}$ here
and $\phi_{R,\epsilon}$
in the previous case is that the term $\frac{\epsilon^l}{|x|^l}$
has been changed to $(\frac1{\log\frac {c|x|}\epsilon})^{\frac1{p-1}}$.
As before, we define
\[
\psi_{R,\epsilon}(x,t)=\phi_{R,\epsilon}(x)\exp(\gamma(t+1)),
\]
and convert  to radial coordinates, with $|x|=r$.
Note that $\frac1{p-1}+1=\frac p{p-1}$ and $\frac1{p-1}+2=\frac{2p-1}{p-1}$.
In place of \eqref{second} and \eqref{other}, we have
\begin{equation}\label{second1}
\begin{aligned}
&\exp(-\gamma(t+1))\big((r-\epsilon)(R-r)\big)^{-\frac2{p-1}-2}
(\psi_{R,\epsilon})_{rr}=\\
&(\frac2{p-1})(\frac2{p-1}+1)(R+\epsilon-2r)^2(1+r)^\frac2{p-1}
\big(1+(\frac{R^2}{\log\frac {cr}\epsilon})^\frac1{p-1}\big)\\
&+2(\frac2{p-1})(r-\epsilon)
(R-r)(1+r)^{\frac2{p-1}}\big(1+(\frac{R^2}{\log\frac {cr}\epsilon})^\frac1{p-1}\big)\\
&-2(\frac2{p-1})^2(r-\epsilon)(R-r)(R+\epsilon-2r)(1+r)^{\frac2{p-1}-1}
\big(1+(\frac{R^2}{\log\frac {cr}\epsilon})^\frac1{p-1}\big)\\
&+(\frac2{p-1})^2(r-\epsilon)(R-r)(R+\epsilon-2r)(1+r)^{\frac2{p-1}}
\frac1{r(\log\frac{cr}{\epsilon})^{\frac p{p-1}}}R^{\frac2{p-1}}\\
&+(\frac2{p-1})(\frac2{p-1}-1)((r-\epsilon)(R-r))^2
(1+r)^{\frac2{p-1}-2}
\big(1+(\frac{R^2}{\log\frac {cr}\epsilon})^\frac1{p-1}\big)\\
&-(\frac2{p-1})^2((r-\epsilon)(R-r))^2
(1+r)^{\frac2{p-1}-1}
\frac1{r(\log\frac{cr}\epsilon)^{\frac p{p-1}}}R^{\frac2{p-1}}\\
&+(\frac 1{p-1})((r-\epsilon)(R-r))^2(1+r)^\frac2{p-1}
\big(\frac1{r^2(\log\frac{cr}\epsilon)^{\frac p{p-1}}}+\frac p
{(p-1)r^2(\log\frac{cr}\epsilon)^{\frac{2p-1}{p-1}}}\big)
R^\frac2{p-1}
\end{aligned}
\end{equation}
and
\begin{equation}\label{other1}
\begin{aligned}
&\exp(-\gamma(t+1))\big((r-\epsilon)(R-r)\big)^{-\frac2{p-1}-2}
\Big(\frac {n-1}r(\psi_{R,\epsilon})_r-\psi_{R,\epsilon}^p
-(\psi_{R,\epsilon})_t\Big)=\\
&-(\frac2{p-1})(\frac {n-1}r)(r-\epsilon)(R-r)(R+\epsilon-2r)
(1+r)^{\frac2{p-1}}\big(1+(\frac{R^2}{\log\frac {cr}\epsilon})^\frac1{p-1}\big)\\
&+(\frac2{p-1})(\frac {n-1}r)((r-\epsilon)(R-r))^2
(1+r)^{\frac2{p-1}-1}\big(1+(\frac{R^2}{\log\frac {cr}\epsilon})^\frac1{p-1}\big)\\
&-(\frac1{p-1})(\frac {n-1}r)((r-\epsilon)(R-r))^2
(1+r)^{\frac2{p-1}}\frac1{r(\log\frac{cr}{\epsilon})^{\frac p{p-1}}}R^{\frac2{p-1}}\\
&-\gamma((r-\epsilon)(R-r))^2(1+r)^{\frac2{p-1}}
\big(1+(\frac{R^2}{\log\frac {cr}\epsilon})^\frac1{p-1}\big)\\
&-(1+r)^{\frac{2p}{p-1}}\big(1+(\frac{R^2}{\log\frac {cr}\epsilon})^\frac1{p-1}\big)^p
\exp((p-1)\gamma(t+1)).
\end{aligned}
\end{equation}

As before, we denote the terms in \eqref{second1} and \eqref{other1}
 by $J_1-J_7$ and $I_1-I_5$
respectively.
It's easy to see that the analysis in the previous case carries over to the present
case when $r$ satisfies $r\ge \delta_0$, where, as above,
$\delta_0$ is an arbitrary positive constant.
It remains to consider $r\in[\epsilon, \delta_0]$.

Exactly as  in the previous case, we have $J_5\le|I_4|$
and $J_2\le|I_4|+|I_5|$,
 and similar to the
previous case, it is easy to see that if $c$ is chosen sufficiently large,
then $J_7\le |I_3|$. (For this last inequality, we use the fact
that the condition $n=\frac{2p}{p-1}$
guarantees that $n>2$.)
We now consider the term $I_2$.
Using the fact that $n-1=\frac{p+1}{p-1}$, and replacing $\frac{r-\epsilon}r$
by 1, we
have
\[
I_2\le\frac{2(p+1)}{(p-1)^2}(r-\epsilon)(R-r)^2
(1+r)^{\frac2{p-1}-1}\big(1+(\frac{R^2}{\log\frac {cr}\epsilon})^\frac1{p-1}\big),
\]
whereas
\[
|J_3|=2(\frac2{p-1})^2\frac{R+\epsilon-2r}{R-r}
(r-\epsilon)(R-r)^2(1+r)^{\frac2{p-1}-1}
\big(1+(\frac{R^2}{\log\frac {cr}\epsilon})^\frac1{p-1}\big).
\]
Since
 $\frac{R+\epsilon-2r}{R-r}$
can be made arbitrarily close to 1 by choosing $R$ sufficiently large,
we have
$I_2\le|J_3|$. (Notice that this argument does not work
in the
case that $n>\frac{2p}{p-1}$ if
$n$ is chosen sufficiently  large. On the other hand, the method
of dealing with $I_2$ that was used above in the case
$n>\frac{2p}{p-1}$---namely, treating it together with $I_1$---does not
work in the present  case that $n=\frac{2p}{p-1}$. It is because
of this that it has been necessary to  split the proof into two cases.)

Now consider the term $J_4$.
In the case that $n>\frac{2p}{p-1}$, $J_4$ was treated together with $I_1$;
in the present borderline case, this will not work. It is here that
the amended form of $\phi_{R,\epsilon}$ is needed.
We have
\[
J_4\le CR^\frac{2p}{p-1}(\log\frac{cr}\epsilon)^{-\frac p{p-1}},
\ \text{for}\ r\in[\epsilon,\delta_0],
\]
for some $C>0$.
On the other hand,
\[
|I_5|\ge MR^{\frac{2p}{p-1}}(\log\frac{cr}\epsilon)^{-\frac p{p-1}},
\ \text{for} \ r\in[\epsilon,\delta_0],
\]
where $M$ can be chosen as large as one wants by choosing
$\gamma$ sufficiently large.
Thus, by choosing $\gamma$ sufficiently large, we have $J_4\le |I_5|$.

Finally, the term $J_1$ is treated as it was in the previous
case, but without the addition of $J_4$ and $I_2$.
Using the fact that $n=\frac{2p}{p-1}$,
the analysis in \eqref{J_1I_1}-\eqref{estimate}
gives
\begin{equation}\label{JI}
J_1+I_1\le
C\frac\epsilon r R(1+r)^{\frac2{p-1}}
(R+\epsilon-2r)\big(1+(\frac{R^2}{\log\frac {cr}\epsilon})^\frac1{p-1}\big).
\end{equation}
Comparing the right hand side of \eqref{JI}
with $|I_5|$, one sees that the inequality
\[
J_1+I_1+I_5\le0, \ \text{for}\ r\in[\epsilon,\delta_0],
\]
will hold with $\gamma$ chosen sufficiently large if
\begin{equation}\label{that'sit}
\frac\epsilon r R^2\le
M\big(1+(\frac{R^2}{\log\frac {cr}\epsilon})^\frac1{p-1}\big)^{p-1},
\ \text{for}\ r\in[\epsilon,\delta_0],
\end{equation}
holds with $M$ chosen sufficiently large.
The right hand side of \eqref{that'sit} is larger than
$MR^2(\log\frac {cr}\epsilon)^{-1}$;
thus, \eqref{that'sit}  holds since
$\frac\epsilon r(\log\frac {cr}\epsilon)$
is bounded for $r\in[\epsilon, \delta_0]$, uniformly over
small $\epsilon$.
This completes the proof of part (2).
\end{proof}

\end{document}